\magnification=1200
\def\dfl{\partial_{\mu}\mkern -18mu\raise 8pt\hbox{$\leftrightarrow$}}

\catcode `@=11
\font\eightrm=cmr8

\font\st=cmr12
\font\fiverm=cmr5
\font\tenbb=msbm10
\newfam\bbfam
\textfont\bbfam=\tenbb

\def\build#1_#2^#3{\mathrel{\mathop{\kern 0pt#1}\limits_{#2}^{#3}}}
\hsize=12cm
\vsize=19cm
\def\tr{\nabla}

\def\overlay#1#2{\ifmmode%
    \setbox0=\hbox{$#1$}%
    \setbox1=\hbox to\wd0{\hss$#2$\hss}\else%
    \setbox0=\hbox{#1}%
    \setbox1=\hbox to\wd0{\hss#2\hss}\fi%
    #1\hskip-\wd0\box1 }

\def\H{\overlay{\bigcirc}{\hbox{\fiverm H}}}

\nopagenumbers

\centerline{ CENTRE DE PHYSIQUE THEORIQUE}
\centerline{ CNRS - Luminy, Case 907}
\centerline{ 13288 Marseille Cedex}
\vskip 3truecm

\centerline{\bf RADIATION STRUCTURES ON ISOTROPIC}
\centerline{\bf HYPERSURFACES}

\bigskip

\centerline{ G. BURDET and M. PERRIN}

\vskip 4truecm
\centerline{\bf Abstract} 

\medskip

The geometrical structures (in the sense of E. Cartan) are analyzed which underlie the gravitational radiation
phenomenon. Among the results are :

- the introduction of the adapted frame bundle to a congruence of isotropic hypersurfaces in a
Lorentzian manifold,

- the description of the reduced frame bundle which admits a unique {\it radiation}
connection induced from the ambiant space-time one, 

- the determination of the automorphisms of an integrable radiation structure,
 
- the repercussions of the geometry on the shape of the stress-energy tensor in the
Einstein's field equations. 

\vfill\eject
\noindent {\st{\bf 1. Introduction}}

\medskip 

In all the standard books of differential geometry the theory
of Riemannian submanifolds of a Riemannian manifold benefits
by a large display having been a field of investigation as old
as the differential geometry itself. In what concerns the pseudo-Riemannian
case authors are usually more discreet and the isotropic degenerate case is a priori
discarded in general. However in spite (or in reason) of their pathological nature,
isotropic submanifolds often appear in physics of the Lorentzian space-time. At least
one can think to the following problems:

\medskip

i) It was early recognized (Mach, Dirac 1938 [1a]) that every cosmological model
should be formulated from data on the past light-cone of the observers, which has led
for instance Dirac [1b] to propose the front form of relativistic dynamics, but also to
study causality and Cauchy development problems using a time parameter whose level
surfaces are isotropic hypersurfaces in space-time [2]. By the way let us note that the tools to
analyse the formation of singularities or the (semi) global existence of solutions of Einstein's
equations are not available yet (for a recent reference see [3]).

In the seventies the same ideas underlie the attempts of particles interaction
modelisations through the infinite momentum frame, the light-cone quantization
and the extreme parton model techniques. 

\medskip

ii) The system of Einstein's equations is of hyperbolic type and as such its characteristics
are hypersurfaces accross which solutions might suffer discontinuities of their derivatives. Then
the radiation concept is formalized as these characteristic isotropic hypersurfaces and the
bicharacteristic rays along which disturbances are propagated. Physically these properties have
led to the gravitational wave concept in vogue during the sixties but next forsaked because of the
non-convincing results of the miscellaneous attempts of detection which followed Weber's
experience. But new elements in favour of existence of gravitational waves were deduced from the
observations of the binary pulsar in 1974 [4], generating a revival of interest for this subject,
and, at the present time, several gravitational wave detection experiments are in preparation as
well as in U.S.A. (cf. the ''Laser Interferometer Gravitational Wave Observatory`` or LIGO
project) as in several European countries (cf the italo-french project VIRGO).

\medskip

iii) In the sixties the fascinating aspect of unexpected features in Einstein's theory came
to development namely the existence of singularities in solutions. An ''outside`` observer does
not perceive the singularity itself but an event horizon. Conversely light cannot escape the
singularity, it becomes infinitely redshifted as the horizon is approached. Then, modulo some
reservations it has been shown that horizons are isotropic hypersurfaces [5].

Hence theoretical
physicists have widely used isotropic hypersurfaces with their degenerate metric induced by the
embedding into the Lorentzian space-time. Now if it is clear that on a manifold with a degenerate
metric there is not a unique related affine connection, concerning isotropic hypersurfaces, due to
the presence of the ambiant Levi-Civit\`a connection, the opinions of authors are divergent (see for
instance [6]).

In the fifties, by taking back the marvellous geometrical Cartan's ideas and essentially under the
impulsion of Erhesman [7] and Chern [8] the notion of connection in fibre bundles and the theory
of $G$-structures have been developed. But these techniques have not been used in the previous
referred theoretical physics works and  consequently the purpose of this paper is to show that a
clear description of connections associated to isotropic hypersurfaces can be given by using the
just above mentioned techniques. Here we do not want to take into account for the presence of
caustics and focal points and we shall restrict ourselves to isotropic hypersurfaces which are
codimension-one differentiable submanifolds of space-time. The paper is organized as follows:

\medskip

\item{-} in Sect.2 the bundle of adapted frames over an isotropic hypersurface is described.
\item{-} In Sect.3 the non-unicity in the bundle of adapted frames of the connection induced by
the Levi-Civit\`a connection of the ambiant space-time is pointed out and it is shown that a unique
induced radiation connection can be defined over a particular subbundle of the bundle of adapted
frames : the so defined radiation structure.
\item{-} Sect.4 is devoted to the descriptions of generalized radiation connections in radiation
structures.
\item{-} In Sect.5 the (infinitesimal) automorphisms of the above introduced radiation structures
are described.
\item{-} A closer inspection of the so-introduced geometrical structures reveals important
constraints on the physical content of the right hand side of the Einstein's field equations. This
point is commented in the conclusion.

\bigskip

\noindent{\st{\bf 2. Bundle of Adapted Frames over a congruence of isotropic hypersurfaces}}

\medskip

Let $(V_{n,1},g)$ denote the Lorentz space-time i.e. a $(n+1)$ - dimensional smooth manifold
endowed with an indefinite metric tensor $g$ with signature $n-1$. In $(V_{n,1},g)$ we want to
consider $n$-dimensional isotropic hypersurfaces $V_n$ with a one-dimensional foliation induced by
an isotropic vector field $\xi$, each leave of which being tangent to an isotropic geodesic : a
bicharacteristic ray. In fact the direction of $\xi$ corresponds to the characteristic rays of the
isotropic hypersurface $V_n$. Therefore, instead of the vector field $\xi$, we have to consider
the line field $[\xi]$ which is the span of the vector field $\xi$ i.e.
$[\xi]=\left\{\lambda\xi,\lambda\in\dot{\bf R}\right\}$, so we have just to suppose that $\xi$ is
both isotropic $(g(\xi,\xi)=0)$ and tangent to a congruence of unparametrized geodesics.

At each point $x$ of $V_n$ the cone of isotropic directions of $x$ has a first
order contact with $V_n$, but the tangent space $T_xV_n$ does not contain any time-like vector,
the future light cone of $x$ being entirely on one side of $V_n$, and the past light cone of $x$
being entirely on the other side. Moreover as an isotropic hypersurface of the Lorentzian ambiant
space-time, $V_n$ inherits of a "degenerate" metric $\beta$, the kernel of which is generated
by the line field $[\xi]$.

Let $G\ell (V_{n,1})$ be the principal fibre bundle of linear frames on $V_{n,1}$. The presence of
the symmetric metric tensor $g$ leads to the reduction of $G\ell (V_{n,1})$ to the bundle of
orthonormal frames $O(V_{n,1})$, the structural group of which being $O(n,1):= \left\{ a\in
G\ell (n+1,{\bf R})\vert ^taSa=S\right\}$, where $t$ denotes the transposition between rows and
columns, and $S$ is chosen as the following $(n+1)\times (n+1)$ symmetric matrix

$$S=\left (\matrix{0&\bf 0&-1\cr\cr
\bf 0&{1\!\!1}_{n-1}&\bf 0\cr\cr
-1&\bf 0&0\cr}\right )\ .\eqno (2.1)$$

 Moreover let us take into account for the presence of the given isotropic line field $[\xi]$ a
representative element of which being written as $^t(0\ { \bf 0 }\  1)$, and look for the subgroup
$G$ of $O(n,1)$ which keeps the isotropic direction fixed. It is easy to see that $G$ is realized by
the matrices of the following form :

$$\left (\matrix{a^{-1}&\bf 0&0\cr\cr
R^t\underline{U}&R&\bf 0\cr\cr
{1\over 2}aU^2&a\underline{U}&a\cr}\right )\eqno (2.2)$$
\noindent where
$R\in O(n-1)$ i.e. $^tRR=1\!\!1_{n-1}$, $a\in \dot{\bf R}$ and
$\underline{U}$ is a $(n-1)$ dimensional row  ($U^2$ denoting the scalar
$^t\underline{U}.\underline{U}$).
\medskip
 Then $G$ can be written as a semi-direct product $(\dot {\bf R} \otimes
O(n-1)) \oslash {\bf R}^{n-1}$. Let us note that the homogeneous space $O(n, 1)/G$ is diffeomorphic to
the $(n-1)$-sphere, and is known as the $(n-1)$-dimensional M$\ddot{\hbox{o}}$bius space in conformal
geometry [9].

Hence one is led to consider the reduction of $O(V_{n,1})$ to a principal fibre bundle with $G$ as
structure group : it will be denoted by $G(V_{n,1})$ and called the frame bundle adapted to the
triplet $(V_{n,1},g,[\xi])$.

Let $(e_0e_1\dots e_{n-1}e_n):= (e_0\ \underline{e}\ e_n)$ denote a moving frame where $e_0$ is
\break isotropic, $\underline{e}$ is a collection of $(n-1)$ space-like vectors and $e_n=\xi$ at
each point of $V_{n,1}$ such that :

$$\matrix{g(e_0,e_0)\hfill&=&g(e_n,e_n)=0\hfill\cr
g(e_A,e_B)\hfill&=&\delta_{AB}\qquad\forall A,B\in [1,n-1]\hfill\cr
g(e_0,e_A)\hfill&=&g(e_n,e_A)=0\hfill\cr
g(e_0,e_n)\hfill&=&-1\hfill\cr}\eqno (2.3)$$

The right action of an element $(a,R,\underline{U})\in G$ on a moving frame is given by :

$$(e_0\ \underline e\  e_n)\longmapsto (a^{-1}e_0+\underline eR\ ^t\underline U + {1\over
2}aU^2e_n\quad\underline eR+a\underline Ue_n\quad ae_n)\eqno (2.4)$$

\noindent and the corresponding action on a dual coframe is written as

$$\theta =\left (\matrix{\theta^0\cr\overline{\theta}\hfill\cr
\theta^n\cr}\right )\longmapsto \left (\matrix{a\theta^0\hfill\cr
-a^t\underline{U}\theta^0+^t\!R\overline{\theta}\hfill\cr
a^{-1}\theta^n-\underline{U}^tR\overline{\theta}+{1\over 2}aU^2\theta^0\hfill\cr}\right ).\eqno
(2.5)$$

Obviously $g=^t\theta S\theta$ is kept invariant under this action.

Now let $ V_n $ be a $ n $-dimensional hypersurface of $ V_{n,1} $ and denote \break $i:V_n\rightarrow
V_{n,1}$ the inclusion map. Then $G(V_{n,1})$ will be called the frame bundle of $V_{n,1}$ adapted to
the hypersurface $V_n$ if, at each point of $ V_n,\  (\underline e \ e_n)$ is a frame of $V_n$.

 Indeed let us consider the bundle $G(V_n)$ induced by $i$ from $G(V_{n,1})$, it is a
principal $G$-bundle over $V_n$ with a homomorphism also denoted by \break $i:G(V_n)\rightarrow
G(V_{n,1})$ which induces $i:V_n\rightarrow
V_{n,1}$ and corresponds to the identity automorphism of $G$. Then, in each point of $V_n$, an
element of $G(V_n)$ is a moving frame of $V_n$. On $V_n$ we have $\theta^0=0$ and the action of
$G$ on a coframe is given by

$$\left (\matrix{\overline{\theta}\hfill\cr
\theta_n\cr}\right )\longmapsto \left (\matrix{^tR\overline{\theta}\hfill\cr
a^{-1}\theta^n-U^tR\overline{\theta}\cr}\right ).\eqno (2.6)$$

 From (2.4) and (2.6) we see that, at each point of $V_n, \ [\xi]=[e_n]$ and $\beta
=^t\overline{\theta}.\overline{\theta}$ are left invariant under the action of $G$. So it is now
clear why $G(V_{n,1})$ has been called the frame bundle adapted to the geometric structure induced
by $i$ over $V_n$ i.e. the structure consisting in the smooth symmetric 2-covariant tensor
field $\beta =i^{\ast}g$ which is degenerate, its kernel being spanned by $\xi$. Now let us show
that $G(V_n)$ is really the right frame bundle to consider over $(V_n,\beta,\xi)$ i.e. that
$G(V_n)$ is a reduction of $G\ell (V_n)$ the bundle of linear frames on $V_n$. First let us
introduce the so-called {\it degenerate orthogonal groups} for $p\leq q$ as follows

$$O^{q-p}(p) : =\left\{ g\in G\ell (q,{\bf R})\vert gS^p(q)^tg=S^p(q)\right\}\eqno (2.7)$$

\noindent where $S^p(q)$ denotes a 2-contravariant symmetric tensor, degenerate of order $q-p$, and

$$O_{q-p}(p) := \left\{ g\in G\ell (q,{\bf R})\vert ^tgS_p(q)g=S_p(q)\right\}\eqno (2.8)$$

\noindent where $S_p(q)$ denotes a 2-covariant symmetric tensor degenerate of order $q-p$. For $p=q$,
$O^0(p)= O_0(p)$ and they are both isomorphic to $O(p)$ the usual orthogonal group.

Let $G\ell(V_n)$ be the bundle of linear frames over $V_n$, it is a principal fibre bundle with
structural group $G\ell$ (dim $V_n, {\bf R })=G\ell ( n,{\bf R})$. Then ${\bf R}^n$ is the standard
fibre of the tangent bundle $T(V_n)$ associated with $G\ell (V_n)$. Since any element
$r\in G\ell (V_n)$ over $x\in V_n$ can be considered as a one-to-one linear mapping of ${\bf R}^n$
onto $T_x(V_n)$ $y\rightarrow ry=Y$, at each point of $V_n$ it is possible to associate with the
metric $\beta$ a bilinear form $(,)_{\beta}$ on ${\bf R}^n$ defined by

$$(y,y')_{\beta}=(r^{-1}Y. r^{-1}Y') =\beta (Y,Y').\eqno (2.9)$$

\noindent This bilinear form can be written

$$(y,y')_{\beta}=^ty\ S_{n-1}(n)y\eqno (2.10)$$

\noindent where $S_{n-1}(n)$ is the $n\times n$ matrix which represents $\beta$ and $y\in{\bf
R}^n$ is written as a column with $n$ elements, $^ty$ denoting the corresponding transposed row.

According to (2.8) the bilinear form $(,)_{\beta}$ is invariant under the action of $O_1(n-1)$.
By choosing

$$S_{n-1}(n)=\left (\matrix{1\!\!1_{n-1}&\bf 0\cr\bf 0&0\cr}\right )\eqno (2.11)$$

\noindent it is easy to verify that $O_1(n-1)$ is a semi-direct product
$\displaystyle (O(n-1)\otimes\dot{\bf R})\oslash{\bf R}^{n-1}$ in which we recognize the
group $G$ previously realized as a subgroup of $G\ell (n+1,{\bf R})$ by the matrices (2.3). Here
$G$ is realized as a subgroup of $G\ell (n,{\bf R})$ by matrices of the form $\displaystyle \left
(\matrix{R&\bf 0\cr a\underline U&a\cr}\right )$.

The invariance of $(,)_{\beta}$ by $O_1(n-1)$ implies that Rel.(2.10) is independent of the choice
of $r$ modulo a right action of an element of $O_1(n-1)$ as a subgroup of $G\ell (n,{\bf R})$ into
$G\ell (V_n)$ i.e. it leads to a reduction of $G\ell (V_n)$ to a $O_1(n-1)$-structure, so we are
led to the following definition :

\medskip

\noindent {\rm Definition} : The bundle of adapted linear frames over $[V_n,\beta ,[\xi ]]$, a congruence
of isotropic hypersurfaces generated by a given line field $[\xi ]$, is a $G$-structure (i.e. a
subbundle $G(V_n)\hookrightarrow G\ell (V_n)$) where $G$ is the degenerate orthogonal group
$O_1(n-1)$.

\medskip

In some particular cases such as the $pp$-waves for instance, a stronger condition is involved,
namely instead of considering $(V_{n,1},g)$ equipped with a line field $[\xi ]$ it must be endowed
with a covariantly constant vector field $\xi$. This leads to introduce another reduction of
$G\ell (V_{n,1})$ to a principal $G_I$-bundle $G_I(V_{n,1})$ where $G_I$ denotes the stabilizer of
$\xi$. According to (2.4) it is deduced from $G$ by setting $a=1$ in (2.3). So the dilation is
excluded and we are left with $G_I=O(n-1)\oslash{\bf R}^{n-1}$.

As previously let us introduce the bundle $G_I(V_n)$ induced from $G_I(V_{n,1})$ by the inclusion
map $i$.
Now $\beta$ and $\xi$ are left invariant under the action of $G_I$ (as a consequence of rel.(2.4)
with $a=1$).

Then from $\xi$ it is possible to construct a degenerate symmetric 2-contrava-riant tensor field
$Z=\xi \otimes \xi$. For the cotangent bundle $T^{\ast}(V_n)$ a construction similar to the one
above described in the case of the tangent bundle can be done. Now at each point $x$ of $V_n$, let
us denote by $Y=ry$ the element of $T^{\ast}_x(V_n)$ corresponding to $y\in{\bf R}^n$. Then an
invariant bilinear form $(,)_{\xi}$ is defined by

$$(y,y')_{\xi}=^ty\ S^1(n)y'=Z(Y,Y')\eqno (2.12)$$

\noindent where $S^1(n)$ is the $n\times n$ matrix which represents $Z$.

Following the previous choice (2.11) let us set 

$$S^1(n)=\left (\matrix{0_{n-1}&\bf 0\cr\bf 0&1\cr}\right ).\eqno (2.13)$$

Then owing to (2.7) it can be shown that the bilinear form $(,)_{\xi}$ is invariant under
$O^{n-1}(1)$ which can be written as $O^{n-1}(1)= G\ell (n-1,{\bf R})\oslash {\bf
R}^{n-1}$.

The group which leaves invariant both bilinear forms $(,)_{\beta}$ and $(,)_{\xi}$ is the
intersection $O_1(n-1)\cap O^{n-1}(1)$ i.e. the group $G_I$ above introduced which can be realized
as a subgroup of $G\ell (n,{\bf R})$ by the matrices of the form

$$\left (\matrix{R&\bf 0\cr \underline U&1\cr}\right ).\eqno (2.14)$$

 So in this case we are led to a $G_I$-structure and to the following definition.

\medskip

\noindent {\rm Definition} : The bundle of adapted linear frames over $[V_n,\beta ,\xi ]$, a
congruence of isotropic hypersurfaces generated by a given vector field $\xi$, is a $G_I$-structure
where $G_I$ is the intersection $O_1(n-1)\cap O^{n-1}(1)$ of two degenerate orthogonal groups.

\bigskip

\noindent {\st{\bf 3. Restriction of the metric connection to the bundle of adapted linear frames
over congruences of isotropic hypersurfaces. The radiation connection}}

\medskip
In the previous section the presence of the metric connection of $(V_{n,1},g)$ has been alluded
through the properties of the vector field $\xi$. Here we want to study the properties of the
metric connection with respect to the bundles of linear adapted frames $G(V_{n,1})$ and
$G_I(V_{n,1})$.

In the bundle of Lorentzian frames $O(V_{n,1})$ let us consider the Levi-Civit\`a connection
$\varphi$ induced by the metric connection. It takes values in the Lie algebra of $O(n,1)$ denoted
 by ${\cal L}(O(n,1))$. For our purpose it is convenient to introduce a dense subset in $O(n,1)$ 
defined as the set of matrices given by the product of three matrix subgroups symbolically
denoted

$$({\bf R}^{n-1})\ (\dot{\bf R}\otimes O(n-1))\ ({\bf R}^{n-1})^{\ast}$$

 \noindent each factor
corresponding to a subgroup of $O(n,1)$  parametrized as follows :

$$({\bf R}^{n-1})\equiv\left\{\left (\matrix{1&^t\overline{V}&{1\over
2}V^2\cr\bf 0&1\!\!1_{n-1}&\overline{V}\cr0&\bf 0&1\cr }\right ),
\overline{V}\  \hbox{column}
\in{\bf R}^{n-1}\right\}$$

$$(\dot{\bf R}\otimes O(n-1))\equiv\left\{\left (\matrix{a^{-1}&&\cr&R&\cr&&a\cr }\right
), \matrix{a\in\dot{\bf R},R\in O(n-1)\hfill\cr ^tRR=1\!\!1_{n-1}\hfill\cr}\right\}$$

$$({\bf R}^{n-1})^{\ast}\equiv\left\{\left (\matrix{1&\bf
0&0\cr^t\underline{U}&1\!\!1_{n-1}&\bf 0\cr {1\over 2}U^2&\underline{U}&1\cr }\right ),
\underline{U}\ \hbox{row
}\in {\bf R}^{n-1}\right\}.$$

Let us note that with $(\dot{\bf R}\otimes O(n-1))\oslash ({\bf R}^{n-1})^{\ast}$ the
parametrization of $G$ used in (2.4) is recovered.

Owing to this symmetric decomposition of a dense subset in $O(n,1)$,  the Lie
algebra ${\cal L}(O(n,1))$ can be decomposed in the following way

$${\cal L}(O(n,1))= {\cal L}({\bf R}^{n-1})+{\cal L}(\dot{\bf R}\otimes O(n-1))+{\cal L}\left
( ({\bf R}^{n-1})^{\ast}\right ) \  , $$

\noindent a vector space direct sum.

Hence the Livi-Civit\`a connection $\varphi$ can be written in the following matrix form

$$\varphi =\left (\matrix{-\phi^n_n&^t\overline{\phi}&0\cr
^t\underline{\phi}&\phi&\overline{\phi}\cr
0&\underline{\phi}&\phi^n_n\cr}\right )\eqno (3.1)$$

\noindent where $\phi^n_n$ is ${\bf R}$-valued, $\overline{\phi}=\left\{\phi^A_n,A\in
[1,n-1]\right\}$ is ${\cal L}({\bf R}^{n-1})$-valued,

\noindent $\underline{\phi}=\left\{\phi^n_A,A\in
[1,n-1]\right\}$ is ${\cal L}\left (({\bf R}^{n-1})^{\ast}\right )$-valued and 
$\phi=\left\{\phi^B_A,A,B\in
[1,n-1]\right\}$ is ${\cal L}\left (O({n-1})\right )$-valued
i.e. such that $\phi^B_A+\phi^A_B=0$.

By introducing the Ricci coefficients, each component of the connection form can be expressed in
terms of the soldering form $\theta$ as follows

$$\varphi^a_b=\gamma^a_{cb}\theta^c\eqno (3.2)$$

\noindent $a,b,c\in [0,n]$, the Ricci's coefficients $\gamma$ being related to the Christoffel's
symbols of the linear connection by

$$\gamma^a_{bc}=e^{\beta}_b\ e^{\gamma}_c\
\left\{\Gamma^{\alpha}_{\beta\gamma}\theta^a_{\alpha}-\partial_{\beta}\theta^a_{\gamma
}\right\}\eqno (3.3)$$

\noindent where $\alpha,\beta,\dots$ are the indices of a local coordinate system.

Then at each point of $V_{n,1}$ the covariant derivative of any vector field $X \ (= X^0 e_0 + X^A e_A
 + X^n e_n )$ is given by :

$$\eqalignno{
\nabla X
= & \ (dX^0  +\ X^A\phi^0_A  +\ X^0\phi^0_0)      \otimes e_0  \cr
  & +\ (dX^C   +\ X^A\phi^C_A   +\ X^0\phi^C_0   +\ X^n\phi^C_n)  \otimes e_C  \cr
  & +\ (dX^n   +\ X^A\phi^n_A                   +\ X^n\phi^n_n)  \otimes e_n&(3.4) \cr} $$

\noindent where $\phi^0_0=-\phi^n_n\ ,\ \phi^0_A=\phi^A_n\ ,\
\phi^A_0=\phi^n_A$.

In the restriction to ${\cal L}(G)$ the components $\phi^0_A$ (hence $\phi^A_n$) disappear. So if
we consider a vector field $X$ tangent to the isotropic hypersurface $V_n$ (i.e. corresponding to
$X^0=0$), rel.(3.4) reduces to

$$\nabla X=(dX^C+X^A\phi^C_A)\otimes e_C+(dX^n+X^A\phi^n_A+X^n\phi^n_n)\otimes e_n\eqno
(3.5)$$

\noindent where there is no component along $e_0$. Therefore for any vector
field $Y$ tangent to $V_n,\ \nabla_YX$ is tangent to $V_n$. But things do not go
on the same way in what concerns one-forms (covariant vectors) because of the
degeneracy of the structure. Indeed, for a one-form $f \  (= f_0 \theta^0 + f_A \theta^A + f_n \theta^n )$,
 the covariant derivative
with respect to the Levi-Civit\`a connection can be written as

$$\eqalignno{
\nabla f
= & \ (df_0  -\ f_A\phi^A_0  -\ f_0\phi^0_0)      \otimes \theta^0  \cr
  & +\ (df_C   -\ f_A\phi_C^A   -\ f_0\phi^0_C   -\ f_n\phi^n_C)  \otimes \theta^C  \cr
  & +\ (df_n-\ f_n\phi^n_n)  \otimes \theta^n\ .&(3.6) \cr} $$

By restriction to a ${\cal L}(G)$-valued connection and for a one-form of \break $V_n \ (f_0=0)$,
one gets

$$\eqalignno{\nabla f=-f_A\phi^A_0\otimes\theta^0&+(df_C-f_A\phi^A_C-f_n\phi^n_C)\otimes
\theta^C\cr
&+(df_n-f_A\phi^A_n-f_n\phi^n_n)\otimes
\theta^n&(3.7)\cr}$$

\noindent in which a $\theta^0$-component remains, showing that $\nabla_Yf\
\forall\ Y\in TV_n$, is not a 1-covariant tensor of $V_n$. Hence the following
proposition is established

\medskip
\noindent{\bf Proposition 3.1}

{\it There is no connection induced on a congruence of isotropic hypersurfaces $[V_n]$ by the
reduction to $G(V_{n,1})$ of the ambiant Levi-Civit\`a connection.}
 \medskip

To make disappear the $\theta_0$ term in (3.7) it is required that $\phi^A_0=0(=\phi^n_A)$. This
leads to restrict the structural group $G$ to a subgroup $G_R$ in which the invariant subgroup
${\bf R}^{n-1}$ has been discarded, that is $G_R=\dot{\bf R}\otimes O(n-1)$. Therefore in a
$G_R(V_{n,1})$-subbundle the covariant derivatives of a contravariant and a covariant vector as
given by (3.5) and (3.7) (with $\phi_0^A=0)$ along a vector field of $V_n$ are
both well defined tensors of $V_n$. So one gets :

\medskip
\noindent {\bf Proposition 3.2}

{\it There is a unique connection induced on $[V_n]$ by the reduction to the bundle $G_R(V_{n,1})$ of
the  Levi-Civit\`a connection of the ambiant space-time. This connection will be called the radiation
connection.} \medskip

A local description of the radiation connection is given in the following section.

But here we want to underline that the stuctural group of any $G$-bundle can be reduced to $G_R$. This
statement follows from the fact that the quotient of $G$ by $G_R$ is contractible, being
topologically equivalent to ${\bf R}^{n-1}$.

\bigskip

\noindent {\st{\bf 4. Radiation connections on degenerate metric structures}}

\medskip
Now we want to discard the ambiant space-time $(V_{n,1},g)$ and to study the degenerate metric
structure $(V_n,\beta,[\xi ])$ per se by keeping in mind the fibre bundles introduced in the
previous section, namely $G(V_n)$ the fibre bundle corresponding to the fibre bundle of adapted
frames on $V_n$, $G_R(V_n)$ which corresponds to the fibre bundle in which there is the radiation
connection and also $G_I(V_n)$ which corresponds to the strict invariance of the vector field $\xi$.

Obviously every connection in $G(V_n)$,  $G_R(V_n)$ or  $G_I(V_n)$ determines a linear
connection of $V_n$, and to keep in memory the origin of these fibre bundles we shall adopt the
following definitions.

\medskip
\noindent{\rm Definition} : A $G$-{\sl radiation} connection is a linear connection with vanishing
torsion induced by a connection in $G(V_n)$.
Alike for $G_R$ and $G_I$-structures respectively.

\medskip
\noindent 4.A. \underbar{$G$-radiation connections.}

\medskip

\noindent{\bf Proposition 4.1}

{\it With respect to any $G$-radiation connection the degenerate metric $\beta$ is parallel.}

\medskip
\noindent{\it Proof} : In a moving frame, $\beta$ can be written as

$$\beta =^t\theta \ {\cal S}\ \theta\eqno (4.1)$$

\noindent where ${\cal S}:=S_{n-1}(n)$ is the $n$-dimensional matrix defined in (2.11).
Then

$$\nabla_c\beta=^t\theta (^t\gamma_c{\cal S}+{\cal S}\gamma_c)\theta\qquad c\in [1,n]\eqno (4.2)$$

\noindent where $\gamma_c=(\gamma^a_{cb})$  denotes the $n\times n$ matrix of
Ricci coefficients introduced in (3.2).

Hence $\nabla_c\beta =0$ is equivalent to $^t\gamma_c{\cal S}+{\cal S}\gamma_c=0$ in which we
recognize, by definition, the Lie algebra ${\cal L}(O_1(n-1))={\cal L}(G)$ . Q.E.D.

\medskip

\noindent{\bf Proposition 4.2}

{\it With respect to any $G$-radiation connection the line field $[\xi ]$ is such that}

$$\tr [\xi]=\chi\otimes [\xi]\eqno (4.3)$$

\noindent {\it where $\chi$ is a real one-form corresponding to the dilation component of the
connection.}

\medskip
\noindent{\it Proof} : The covariant derivative corresponding to a linear connection in a moving
frame verifies

$$\tr_{e_a}e_n=\gamma^b_{an}e_b.\eqno (4.4)$$

\noindent For a $G$-radiation connection $\phi^B_n=0$ so that $\gamma^B_{an}=0\ \forall B\in
[1,n-1]$. Hence

$$\tr_{e_a}e_n=\gamma^n_{an}e_n.\eqno (4.5)$$

\noindent So, in a moving frame in which $\xi = e_n$, it can be set $\gamma^n_{an}=\chi_a$ and
then

$$\tr_a\xi =\chi_a\xi .\eqno (4.6)$$

\noindent{\bf Proposition 4.3}

{\it With respect to a $G$-radiation connection}

\medskip

{\it i) the line field $[\xi ]$ is geodesic} $\tr_{\xi}\xi = \chi(\xi )\xi$
\medskip
{\it ii) the expansion of $[\xi ]$ does not vanish, div }$\xi = \chi(\xi )$.

\medskip

The proof is obvious from (4.5) and let us note that div $\xi =\chi_n$.

\medskip

\noindent {\bf Proposition 4.4}

{\it The structure equations and Bianchi's identities of a $G$-radiation connection are given by}

$$\matrix{\hbox{(a)}\ &0 = \overline{\H} = d\overline{\theta} + \phi\wedge
\overline{\theta}\hfill\cr
 \hbox{(b)}\ &0 = \H^n=d\theta^n+\underline{\phi}\wedge\overline{\theta}+\phi^n_n\wedge\theta^n\cr
\hbox{(c)}\ &\Phi =d\phi +\phi\wedge\phi\hfill \cr
\hbox{(d)}\ &\underline{\Phi} = d\underline{\phi} + \underline{\phi}\wedge (\phi -
1\!\!1_{n-1}\phi^n_n )\hfill\cr
\hbox{(e)}\ &\Phi^n_n = d\phi^n_n\hfill\cr}
\quad\matrix{\hbox{(a')}\ &\Phi\wedge\overline{\theta}=0\hfill\cr
\hbox{(b')}\ &\underline{\Phi}\wedge\overline{\theta}+\Phi^n_n\wedge\theta^n=0\hfill\cr
\hbox{(c')}\ &D\Phi =0\hfill\cr
\hbox{(d')}\ &D\underline{\Phi}=0\hfill\cr
\hbox{(e')}\ & D\Phi^n_n=0\hfill\cr}\eqno (4.7)$$

\medskip

\noindent{\it Proof} : By using the same conventions as in (3.1) any $G$-radiation connection can be
written under the matrix form

$$\left (\matrix{\phi&\bf 0\cr\underline{\phi}&\phi^n_n}\right )\eqno (4.8)$$

\noindent from which the torsion and curvature 2-forms are easily deduced and the Bianchi's identities
follow.
\medskip
If $\{\overline{\theta},\theta_n\}:=\vartheta$ denote the soldering form of $G(V_n)$ every component
of the connection can be decomposed by using Ricci's coefficients as previously (see (3.2)). The
antisymmetry of the ${\cal L}(O(n-1))$-valued component $\phi$ leads to :

$$\gamma^A_{CB}=-\gamma^B_{CA}\qquad \hbox{and}\qquad \gamma^A_{nB}=-\gamma^B_{nA}\eqno (4.9)$$

Let us now study if among all existing $G$-radiation connections one of them is a privilegied one,
in other words if the miracle of the Riemannian geometry occurs again.

\vfill\eject
\noindent {\bf Proposition 4.5}

{\it All $G$-radiation connections have a common orthogonal component.}

\medskip

\noindent{\it Proof} : To compare the ${\cal L}(O(n-1))$-components of two $G$-radiation
connections $\Gamma$ and $\Gamma '$ let us set $\Delta\phi :=\phi '-\phi =(\gamma '-\gamma )\theta
:=(\Delta\gamma )\theta$.

Then from (4.7a) one can deduced

$$\hbox{(a)}\quad \Delta\gamma^A_{[BC]}=0\qquad \hbox{(b)}\quad \Delta\gamma^A_{nB}=0\eqno (4.10)$$

\noindent and from (4.7b)

$$\hbox{(a)}\quad \Delta\gamma^n_{[AB]}=0\qquad 
\hbox{(b)}\quad \Delta\gamma^n_{[nB]}=0\eqno (4.11)$$

Then (4.7a) together with the antisymmetry properties (4.9) leads to $\Delta\gamma^A_{BC}=0\ \forall
A,B,C$. As $\Delta\gamma^A_{nB}=0$ it can be deduced $\Delta\phi^A_B=0\ \forall A $and$ B$ that is to say
 all $G$-radiation connections have the same ${\cal L}(O(n-1))$-component. In what concerns the
other components it remains

$$\left\{\matrix{
\hbox{(a)}&\Delta\phi^n_B=\Delta\gamma^n_{AB}\theta^A\cr
\hbox{(b)}&\Delta\phi^n_n=\Delta\gamma^n_{An}\theta^A+\Delta\gamma^n_{nn}\theta^n\cr}\right. .\eqno
(4.12)$$

\medskip

\noindent{\bf Proposition 4.6}

{\it On the degenerate structure $(V_n,\beta ,[\xi ])$ the $G$-radiation connections are in
correspondence with the symmetric 2-contravariant tensors on $V_n$.}
\medskip

\noindent{\it Proof} : It is a consequence of (4.11).
\medskip
Let us give the corresponding local expressions. In terms of its Ricci's coefficients and in a
local coordinate system the Christoffel's symbols of a $G$-radiation connection are given by

$$\Gamma^{\alpha}_{\beta\gamma}=e^{\alpha}_a\left\{ \partial_{\beta}\theta^a_{\gamma }+
\theta^b_{\beta}\theta^c_{\gamma }\gamma^a_{bc}\right\}.\eqno (4.13)$$

By performing a standard algebra one gets

$$\Gamma^{\alpha}_{\beta\gamma}=(\underline{e}\otimes
\underline{e})^{\alpha\rho}B_{\rho,\beta\gamma}+e^{\alpha}_n\Gamma^{\rho}_{\beta\gamma}
\theta^n_{\rho}\eqno (4.14)$$

\noindent where $B_{\rho,\beta\gamma}$ is recognized as a Koszul's term defined by

$$\eqalignno{B_{\rho,\beta\gamma}&= (\overline{\theta}_{(\beta}.\partial_{\gamma
)}\overline{\theta}_{\rho}) + (\overline{\theta}_{\rho}
.\partial_{(\beta}\overline{\theta}_{\gamma )})-(\overline{\theta}_{(\beta}.\partial_{\rho}
\overline{\theta}_{\gamma )})\cr
&\equiv \partial_{(\beta}\beta_{\rho\gamma )}-{1\over 2}
\partial_{\rho}\beta_{\beta\gamma}.&(4.15)\cr}$$

In the expression (4.14) a two-covariant tensor $\underline e\otimes \underline e$ is arising, it is
degenerate its kernel being spanned by $\theta^n$. This leads to locally introduce a quasi inverse
$\beta_f$ of $\beta$ (i.e. a two-covariant tensor $\beta_f$ such that its contraction with $\beta\
\beta_f$ gives $\beta_f$) associated to the one-form $f$ such that $f(\xi)=1$ and here identified with  $\theta^n$. More  exactly $\beta_f$ is associated to the choice of a projective one-form $[f]$
such that $[f]([\xi ])=1$, in such a way one can write

$$\beta\ \beta_f=1\!\!1 -[f]\otimes [\xi ]. \eqno (4.16)$$

It is then easy to verify that the general expression which verifies (4.14) can be written under
the following form :

$$\Gamma^{\alpha}_{\beta\gamma}=\beta^{\alpha\rho}_f B_{\rho
,\beta\gamma}+\xi^{\alpha}Z_{\beta\gamma}\eqno (4.17)$$

\noindent where $Z_{\beta\gamma}$ is an arbitrary 2-contravariant symmetric tensor (cf.
proposition 4.5).

Furthermore we have the following local relations :

$$\chi_{\alpha}=-\xi^{\lambda}\tr_{\alpha}f_{\lambda}=-\xi^{\lambda}(\partial_{\alpha}f_{\lambda}
-Z_{\alpha\lambda})\eqno
(4.18)$$

$$\tr_{[\gamma}\chi_{\beta ]}={1\over 2}\xi^{\alpha}R^{\lambda}_{\alpha [\beta\gamma
]}f_{\lambda}\eqno (4.19)$$

\noindent where the $R^{\lambda}_{\alpha\beta\gamma}$'s denote the local components of the
Riemannian curvature tensor which can be written

$$R^{\lambda}_{\alpha\beta\gamma}=\build{R}_{}^{o}\!\!\!\
^{\lambda}_{\alpha\beta\gamma}+2\xi^{\lambda}\left\{(\partial_{[\gamma} +\chi_{[\gamma})Z_{\beta
]\alpha}+\build{\Gamma}_{}^{o}\!\!\!\ ^{\sigma}_{\alpha [\beta}Z_{\gamma ]\sigma}\right\}$$

\noindent where $\displaystyle\build{R}_{}^{o}\!\!\!\
^{\lambda}_{\alpha\beta\gamma}$ denotes the
contribution corresponding to the Koszul's term (i.e. which corresponds to $Z=0$ in (4.17)).

\medskip
\noindent 4.B. \underbar{$G_R$-radiation connections.}

As suggested by results of subsection 4.A. let us briefly discuss the linear connections defined
from ${\cal L}(G_R)$-valued connections we shall name $G_R$-radiation connections.

All the calculations of subsection 4.A can be performed again after having set
$\underline{\phi}=0$ i.e. $\gamma^n_{aB}=0$. Then Propositions 4.1-2-3-5 remain unchanged and
structure equations and Bianchi's identities are deduced directly from (4.7). Therefore Prop 4.6 is
modified and  replaced by the following one.

\medskip

\noindent{\bf Proposition 4.7}

{\it On the degenerate structure $(V_n,\beta,[\xi ])$ the $G_R$-radiation connections are in
correspondence with the one-forms on $V_n$.}
\medskip

Indeed from (4.12b) the difference between the dilation components of two $G_R$-radiation
connections is given by $\Delta\phi^n_n=\Delta\gamma^n_{nn}\theta^n$, which shows that they are not
related. Q.E.D.

For the Christoffel's symbols, a direct calculation from (4.13) leads to :

$$\build{\Gamma}_{}^{o}\!\!\!\ ^{\alpha}_{\beta\gamma}=\beta^{\alpha\rho}_fB_{\rho
,\beta\gamma}+\xi^{\alpha}\left\{\partial_{(\beta}f_{\gamma )}+\chi_{(\beta}f_{\gamma
)}\right\}.\eqno (4.20)$$

\noindent  On $(V_n,\beta ,[\xi ])$ the choice of a quasi-inverse $\beta_f$ of
$\beta$ associated with a one-form $f$ does not fully determine a $G_R$-radiation connection, a
supplementary one-form $\chi$ is needed.

Now we can go back to the unique induced radiation connection defined in proposition (3.2). Indeed
the local expression of the radiation connection over an isotropic hypersurface is also given by
(4.20) in which the one-form $\chi$ corresponds to the pull-back of the dilatation component of the
ambiant Levi-Civit\`a connection. While we are on the subject let us mention that the three
references [5a,c,d] deal with particular choices of $G$-radiation connections (for $n=3$) which
therefore cannot be identified with the unique induced radiation connection.

\medskip
\noindent 4.C. \underbar{$G_I$-radiation connections.}

As mentioned in sect.2, to require the invariability of the vector field $\xi$ leads to
restrict the group $G$ to its subgroup $G_I$ and then to consider \break $G_I$-radiation
connections. They always satisfy proposition 4.1 but propositions 4.2 and 4.3 are replaced by the
following one.

\medskip
\noindent{\bf Proposition 4.8}

{\it With respect to any $G_I$-radiation connection $\xi$ is covariantly constant \break $\tr\xi
=0$, hence geodesic $\tr_{\xi}\xi =0$ and divergence free div$\xi =0$.}

\medskip
\noindent{\it Proof} : It is trivially deduced from the proof of proposition 4.2 by setting
$\gamma^n_{an}=0\ \forall a\in [1,n]$ i.e. by accounting for $\phi^n_n=0$.

Structure equations and Bianchi's identities are deduced from (4.7). Pro-position 4.5 is still valid
for $G_I$-radiation connections. Now Rel.(4.11b) being replaced by  $\Delta\gamma^n_{nB}=0$
this leads to rewrite Prop.4.6 as follows.

\medskip
\noindent{\bf Proposition 4.9}

{\it On the degenerate structure $(V_n,\beta , \xi )$ the $G_I$-radiation connections are in
correspondence with the degenerate symmetric 2-contravariant tensors on $V_n$ the kernel of which
is generated by $\xi$.}

\bigskip

\noindent{\st{\bf 5. Automorphisms of radiation structures}}

\medskip
Among the elements of Dif$(V_n)$ we have to select those which preserve the degenerate structure
$(\beta ,[\xi ])$ that is the transformations of $V_n$ which leave $\beta$ and $[\xi ]$ invariant
in the sense that they transform a representative element of $[\xi ]$ into another one. Obviously
such a transformation $f$ maps each adapted frame at an arbitrary point $x\in V_n$ into an adapted
frame at $f(x)\in V_n$, that is to say that there exists an induced transformation $\widetilde{f}$
of $G\ell (V_n)$ which maps $G(V_n)$ into itself.

However due to the absence of a privilegied $G$-radiation connection in one-to-one correspondence
with $(\beta ,[\xi ])$, in order to promote $f$ or $\widetilde{f}$ to the status of radiation
transformation we have to prescribe that they preserve the chosen $G$-radiation
connection (affine transformation).

\medskip
\noindent{\rm Definition} : A transformation $f$ of $V_n$ is a {\sl radiation transformation} if the
induced transformation $\widetilde{f}$ of $G\ell (V_n)$ maps $G(V_n)$ endowed with a chosen
radiation connection into itself.

Such a definition is fully adapted to define a notion of equivalence between two $G_R$
(or $G_I$)-radiation connections.

\medskip

\noindent{\rm Definition} : Two $G_R$ (respectively $G_I$)-structures are said {\sl equivalent} if the
two considered $G_R$ (respectively $G_I$)-radiation connections are subordinate to the same
$G$-connection.

In fact two different $G_R$ (respectively $G_I$)-structures are associated to two distinct
embeddings $h$ and $h'$ of $G_R(V_n)$ (respectively $G_I(V_n))$ into $G(V_n)$. Then on each fibre
of the bundle $G(V_n)$ these two embeddings $h$ and $h'$ define two $G_R$-orbits (respectively
$G_I$-orbits) such that $h'(p)=h(p)\lambda (x)$ for any $p\in G_R(V_n)$ (respectively $G_I(V_n)$)
which projects onto $x\in V_n,\lambda (x)$ being identified with an element of $G$ which can be
written under the form 

$$\left (\matrix{1\!\!1_{n-1}&\bf 0\cr \bf 0&\rho (x)\cr}\right )= \lambda_R(x)\left
(\hbox{respectively}\ \left (\matrix{1\!\!1_{n-1}&\bf 0\cr \underline{\rho}(x)&1\cr}\right
)=\lambda_I (x)\right )$$

\noindent where $\rho$ is a mapping of $V_n$ into the subgroup $\dot{\bf R}\subset G$ and
$\underline{\rho}$ is a mapping of $V_n$ into the subgroup ${\bf R}^{n-1}\subset G$.

Under a radiation transformation $\widetilde f$ a $G_R$-structure (respectively a $G_I$-structure)
is transformed into an equivalent one.

\medskip

\noindent{\rm Definition} : A vector field $X$ on $V_n$ is called an {\sl infinitesimal radiation
transformation} if the local one-parameter group of local transformations generated by $X$ in a
neighborhood of each point of $V_n$ consists of local radiation transformations.

\medskip

\noindent{\bf Proposition 5.1}

{\it With respect to an infinitesimal radiation transformation one has}

\medskip

{\it i) $L_X\beta =0$ where $L_X$ denotes the Lie derivation with respect to $X$ , }
\medskip
{\it ii) $L_X\xi=[X,\xi ]=k\xi$ with $k=$constant , }
\medskip
{\it iii) $K(X,Y)=0$ for all vector fields $Y$ on $V_n$ where}

\noindent $K(X,Y)=R(X,Y)-\tr_YA_X$, $A_X$
{\it denoting the derivation defined by $L_X-\tr_X$, $R$ being the curvature tensor and $\tr$ the
covariant derivation corresponding to the chosen  $G$-radiation connection.}
\medskip
This property is a direct consequence of the above definitions and iii) defines an infinitesimal
affine transformation [9].

Let us denote by $\widetilde X$ the vector field on $G\ell (V_n)$ induced by  the local group of
 radiation transformations prolonged to $G\ell (V_n)$.

\medskip
\noindent{\bf Proposition 5.2}

{\it Corresponding to any infinitesimal radiation transformation $X$ there is an infinitesimal
automorphism $\widetilde X$ of $G(V_n)$.}
\medskip

This proposition is a direct consequence of the definition of an infinitesimal radiation
transformation.

\medskip
\noindent{\bf Proposition 5.3}

{\it The Lie derivatives with respect to $\widetilde X$ of the canonical one-form 
$\vartheta =\{\overline{\theta},\theta^n\}$ of $G\ell (V_n)$ reduced to a
$G$-structure and those of the chosen $G$-radiation connection one-form
$\varphi
=\left\{\phi ,\phi^n_n ,\underline{\phi}\right\}$  satisfy the following properties}

$$\matrix{L_{\widetilde X}\varphi =0&&\cr L_{\widetilde X}\overline{\theta} =0&,&L_{\widetilde
X}\theta^n =k\theta^n\cr}$$

\noindent {\it where $k$ is a constant.}
\medskip

Finally by comparing the right action of $G$ with those of $G_R$ ($G_I$ respectively) the
following properties can be deduced.

\medskip

\noindent{\bf Proposition 5.4}

{\it The Lie derivatives with respect to $\widetilde X$ of the one-form $\varphi_0=\{\phi
,\phi^n_n\}$ of a chosen $G_R$-radiation connection and of the canonical one-form $\vartheta$ of
$G\ell (V_n)$ reduced to a $G_R$-structure satisfy}

$$\matrix{L_{\widetilde X}\phi =0&,&L_{\widetilde X}\phi^n_n =0\cr L_{\widetilde
X}\overline{\theta}=0&,&L_{\widetilde X}\theta^n=-\underline{\epsilon}\overline{\theta}\cr}$$

\noindent {\it where $\underline{\epsilon}$ is the infinitesimal gauge transformation corresponding
to the above introduced mapping $\underline{\rho}$.}

\medskip
\noindent{\bf Proposition 5.5}

{\it The Lie derivatives with respect to $\widetilde X$ of the one-form $\varphi_1=\{\phi
,\underline{\phi}\}$ of a chosen $G_I$-radiation connection and of the canonical one-form
$\vartheta$ of $G\ell (V_n)$ reduced to a $G_I$-structure satisfy}

$$\matrix{L_{\widetilde X}\phi =0&,&L_{\widetilde X}\underline{\phi}=-\epsilon\underline{\phi}\cr
L_{\widetilde X}\overline{\theta}=0&,&L_{\widetilde
X}\theta^n=-\epsilon\theta^n\cr}$$

\noindent {\it where $\epsilon$ is the infinitesimal gauge transformation corresponding to the
mapping $\rho$.}
\medskip

To illustrate the above established propositions let us treat the case of the standard radiation
structure on ${\bf R}^n$ in which $G(V_n)$ is diffeomorphic to the trivial bundle ${\bf R}^n\times
G$ endowed with the natural $G$-connection provided by the Maurer-Cartan form of $G$. The Lie
algebra ${\cal L}(G)$ clearly contains an element of rank one and hence is of infinite type. So
the general theorem on the automorphisms group of a $G$-structure cannot be applied and it cannot
be claimed that it is a Lie group. Then let us perform the calculation of the infinitesimal
automorphisms of this structure which firstly amounts to select among the
vector fields $X$ of $V_n$ those which satisfy i) and ii) of Prop. (5.1). In a special admissible
coordinate frame system 

$$(x^1,\dots ,x^n):=(\overline x,x^n)\equiv (\{x^A\},x^n(A\in [1,n-1]))$$
\noindent in which $\theta^A_{\alpha}=\delta^A_{\alpha}$ and $e^{\gamma}_n=\delta^{\gamma}_n$ so
that $\beta_{\alpha\gamma}=\theta^A_{\alpha}\theta^C_{\gamma}\delta_{AC}$ and
$\xi^{\gamma}=\delta^{\gamma}_n$, the equations i) and ii) of Prop. (5.1) can be written as

$$\matrix{L_X\beta_{\alpha\gamma}&=&2\tr_{(\alpha}X_{\gamma )}&=&0\cr
L_X\xi^{\gamma}&=&-\xi^{\alpha}\partial_{\alpha}X^{\gamma}&=&k\xi^{\gamma}\ .\cr}$$

\noindent Their general solution is given by :

$$\left\{\matrix{X^A&=&\omega^A_B(\overline x )x^B+a^A\cr
X^n&=&-kx^n+f(\overline x )\cr}\right.\eqno (5.1)$$

\noindent where $\omega^A_B(\overline x )=-\omega^B_A(\overline x)$, $\{a^A\}$ is a set of $n-1$
constants and $f(\overline x)$ can be any function of $\overline x$ so that the corresponding Lie
algebra is infinite dimensional. But by taking iii) of the Prop. (5.1) into account we are led to
select linear (in $\overline{x}$) vector fields only.

Hence the vector fields satisfying i),ii) and iii) of Prop. (5.1) can be written as

$$\left\{\matrix{X^A&=&\omega^A_Bx^B+a^A\cr
X^n&=&k_Bx^B+kx^n+k_0\cr}\right.\eqno (5.2)$$

The Lie brackets of these vector fields generate a Lie algebra which is \break recognized as
the Lie algebra of the inhomogeneous $G$ group isomorphic to ${\bf R}^n\oslash\left (
(O(n-1)\otimes \dot{\bf R})\oslash{\bf R}^{n-1}\right )$.

Now if we recall that a $G$-structure $G(V_n)$ is said integrable if it is locally isomorphic to
the standard $G$-structure on ${\bf R}^n$, from the above result the following proposition can be
 deduced immediately.

\medskip
\noindent{\bf Proposition 5.6} 

{\it The infinitesimal automorphisms of an integrable radiation structure\break $G(V_n)$ endowed
with a radiation connection generate a Lie algebra of dimension $\displaystyle {\eightrm 1\over
2}(n^2+n+2)$ isomorphic to the Lie algebra of the inhomogeneous $G$ group namely}
${\cal L}\left ({\bf R}^n\oslash\left (
(O(n-1)\otimes\dot{\bf R})\oslash{\bf R}^{n-1}\right )\right )$.
\medskip

As a final remark let us compare the automorphisms of a degenerate structure and those of an
isotropic hypersurface. Among the isometry transformations of $(V_{n,1},g)$ we have to select the
ones $F$ which do not move the isotropic rays generated by $[\xi ]$ (or $\xi$ under the covariant
constancy hypothesis). Then Aut $G(V_{n,1})$ is given by the induced transformations $\{\widetilde
F\}$ of $G\ell (V_{n,1})$ which map $G(V_{n,1})$ onto itself. Let us recall that $G(V_n)$ can be
obtained as the pull-back of $G(V_{n,1})$ induced by the inclusion map $i:V_n\rightarrow V_{n,1}$.
Hence to obtain Aut$(V_n)$ we have to select the subgroup of the automorphisms $\{\widetilde F_0\}$
induced by the diffeomorphisms $\{F_0\}$ of $V_{n,1}$ which keep invariant the isotropic
hypersurfaces.

In the above described example of a standard radiation structure embedded into the Minkowski
space-time ${\bf R}^{n,1}$, only one translation which maps an isotropic hyperplane into another
one has to be removed from the subgroup

$${\bf R}^{n+1}\oslash\left ( (O(n-1)\times\dot{\bf R})\oslash{\bf
R}^{n-1}\right )\subset {\bf R}^{n+1}\oslash O(n,1)$$

\noindent Then the group of automorphisms of an integrable radiation structure is recovered.

\medskip
\noindent {\st{\bf 6. Conclusion}}

In this section we restrict ourselves to the case $n=3$. We have seen so far that the existence of
a congruence of isotropic hypersurfaces in a Lorentzian space-time $(V_{3,1},g)$ implies the
reduction of the orthogonal frame bundle\break $O(V_{3,1})$ to the bundle $G(V_{3,1})$ of
adapted frames. We have then to handle with $G$-connections, the curvature forms of which are ${\cal
L}(G)$-valued. But they have also to satisfy the Bianchi's identities leading to the disappearance
of some components in the Ricci tensor, only seven of them are precisely different from zero.
Consequently the physical right hand side tensor $T$ in the Einstein's field equations cannot be
anything. It is easy to verify that for a $G$-connection the compatible physical tensor should be
written in an adapted frame as 

$$\eqalignno{T_G=&\lambda \left (\theta^1\otimes\theta^1+\theta^2\otimes\theta^2
\right)+\pi\theta^0\otimes
\theta^0+\sigma\left(\theta^1\otimes\theta^2+\theta^2\otimes\theta^1\right)\cr
&+ \rho_1(\theta^1\otimes\theta^0+\theta^0\otimes\theta^1)+
\rho_2(\theta^2\otimes\theta^0+\theta^0\otimes\theta^2)+
\rho_3(\theta^3\otimes\theta^0+\theta^0\otimes\theta^3)\cr}$$

\noindent where $\lambda , \pi , \sigma , \rho$ are arbitrary functions on 
$V_{3,1}$ a priori. We do not try to give an interpretation of $T_G$ as a whole
and directly perform the reduction to $G_I(V_{3,1})$ which corresponds to keep
the isotropic vector field $\xi$ covariantly constant. The disappearance of the
${\bf R}^{n-1}$-subalgebra leads to $R_{03}=R_{30}=0$. Consequently the
corresponding stress-energy tensor $T_{G_I}$ is deduced from the above
expression of $T_G$ by setting $\lambda = 0=\sigma$ and $\rho_3 = {1\over 2} 
{\bf S},{\bf S}$ denoting the scalar curvature. Then $T_{G_I}$ can be interpreted as
the stress-energy tensor [10] of a massless particle beam with possible heat flow
$\{\rho_1 ,\rho_2,\rho_3\}$ along the isotropic hypersurface, the radiation
phenomena.

Finally in the case of the reduction to $G_R(V_{3,1})$ which makes appear the uniquely induced
radiation connection, the Ricci tensor is skinny, two components being
non-vanishing only. Then the expression of $T_{G_R}$ involves the components
$\lambda$ and $\rho_3$ related by $\rho_3-\lambda={1\over 2} 
{\bf S}$, and must be interpreted as the energy
tensor of the vacuum [11][12].

 \vfill\eject
 \noindent {\bf REFERENCES}
\medskip
\parindent=1cm

\item{\hbox to\parindent{\enskip [1]}\hfill}P.A.M. Dirac, a) Proc. Roy. Soc. Lond. {\bf A165}, 199-208 (1938);
\item{\hbox to\parindent{\enskip }\hfill}\hphantom{P.A.M. Dirac,} b) Rev. Mod. Phys. {\bf 26}, 392-399 (1949).

\item{\hbox to\parindent{\enskip [2]}\hfill}R. Penrose, {\it Null Hypersurface Initial Data for
Classical Fields of Arbitrary Spin and for General Relativity}, Published as Golden
Oldie in Gen. Rel. Grav. {\bf 12}, 225-264  (1980);

\item{\hbox to\parindent{\enskip }\hfill}H. Bondi, M. Vanderburg, A. Metzner, {\it Gravitational
Waves in General Relativity. VII. Waves from Axi-Symmetric Isolated Systems}, Proc. Roy. Soc. Lond.
{\bf A269},  21 (1962);

\item{\hbox to\parindent{\enskip }\hfill}R.K. Sachs, {\it Gravitational Waves in General Relativity.
VIII. Waves in Asymptotically Flat Space-Time}, Proc. Roy. Soc. Lond. {\bf A270},  103-126 (1962).

\item{\hbox to\parindent{\enskip [3]}\hfill} P.T. Chrusciel, {\it Semi-Global Existence and
Convergence of Solutions of the Robinson-Trautman (2-dimensional Calabi) Equation},
Commun.\break Math. Phys. {\bf 137}, 289-313 (1991).

\item{\hbox to\parindent{\enskip [4]}\hfill}J.H. Taylor and P.M. McCulloch, Ann. of the N.Y.A.S.
{\bf 336}, 442-446  (1980).

\item{\hbox to\parindent{\enskip [5]}\hfill}S.W. Hawking, {\it The Event Horizon}, Lecture at the
Summer School, Les Houches, C. Dewitt and B.S. Dewitt Ed. (1972).

\item{\hbox to\parindent{\enskip [6]}\hfill}a) G. Lemmer, {\it On Covariant Differentiation Within a
Null Hypersurface}, Il Nuovo Cimento {\bf 37}, 1659-1672  (1965) ;

\item{\hbox to\parindent{\enskip }\hfill}b) J.B. Kammerer, {\it Tenseur de courbure d'une
hypersurface isotrope}, C.R. Acad. Sc. Paris {\bf 264}, 86-89  (1967);

\item{\hbox to\parindent{\enskip }\hfill}c) G. Dautcourt, {\it Charateristic Hypersurfaces in
General Relativity}, Jour. Math. Phys. {\bf 8}, 1492-1501  (1967);

\item{\hbox to\parindent{\enskip }\hfill}d) P. Hajicek, {\it Exact Models of Charged Black Holes. 1.
Geometry of Totally Geodesic Null Hypersurface}, Commun. Math. Phys. {\bf 34}, 37-52  (1973).

\item{\hbox to\parindent{\enskip }\hfill}e) A. Ashtekar, M. Steubel, {\it Symplectic Geometry of
Radiative Modes and Conserved Quantities at Null Infinity}, Proc. Roy. Soc. Lond. {\bf A376},
585-607  (1981).

\item{\hbox to\parindent{\enskip [7]}\hfill}C. Ehresmann, {\it Les connexions infinit\'esimales dans
un espace fibr\'e\break diff\'erentiable}, in Colloque de Topologie, Bruxelles, pp. 29-55, Thone,
Li\`ege (1950).

\item{\hbox to\parindent{\enskip [8]}\hfill}S.S. Chern, {\it G-Structures}, in Colloquium, Strasbourg,
pp. 119-136,\break CNRS Paris (1953).

\item{\hbox to\parindent{\enskip [9]}\hfill}S. Kobayashi, Transformation Groups in Differential
Geometry, p.133, \break Springer Verlag Ed. (1972).

\item{\hbox to\parindent{\enskip [10]}\hfill}C.W. Misner, K.S. Thorne, J.A. Wheeler, Gravitation,
Ed. Freeman , San Francisco (1973) ;

\item{\hbox to\parindent{\enskip }\hfill}R.K. Sachs, H. Wu, General Relativity for Mathematicians,
Springer-Verlag (1977).

\item{\hbox to\parindent{\enskip [11]}\hfill}Ya.B. Zel'dovich, I.D. Novikov, The Structure and
Evolution of the Universe, University of Chicago Press (1983).
 
\item{\hbox to\parindent{\enskip [12]}\hfill}G. Burdet, M. Perrin, {\it  Gravitational waves without 
gravitons}.	Letters in Mathematical Physics {\bf 25}, 39-45 (1992).

\end